\documentclass[11pt,a4]{amsart}
\usepackage{amscd,amsmath,amssymb,amsfonts,ifthen,hyperref,tikz-cd} 
\usepackage{calligra}
\usepackage{mathtools}

\usepackage{bbold,tikz}
\usepackage{color} 
\usepackage{graphicx}
\usepackage[all]{xypic}

\usepackage{fourier}

\usepackage[T1]{fontenc}
\newtheorem{claim}{Claim}[section]
\newtheorem{thm}[claim]{Theorem}
\newtheorem{lem}[claim]{Lemma}
\newtheorem{cor}[claim]{Corollary}

\newtheorem{prop}[claim]{Proposition}
\theoremstyle{definition}

\newtheorem{rem}[claim]{Remark}

\usepackage{geometry}
 \date{June,  2024}
 
 \author[Bao]{V\~o Qu\^oc Bao}
 \email[V\~o Qu\^oc Bao]{vqbao@math.ac.vn}
 \author[Hai]{Ph\`ung H\^o Hai}
 \email[Ph\`ung H\^o Hai]{phung@math.ac.vn}
 \author[Thinh]{D\`ao Van Thinh}
 \email[D\`ao Van Thinh]{dvthinh@math.ac.vn}
 \address[V\~o Qu\^oc Bao, Ph\`ung H\^o Hai, D\`ao Van Thinh]{Institute of Mathematics, Vietnam Academy of Science and Technology.}  

\title[Cohomology of the differential fundamental group]{Cohomology of the differential fundamental group of algebraic curves}

\makeatletter
\@namedef{subjclassname@2020}{\textup{2020} Mathematics Subject Classification}
\makeatother

\thanks{The work of V\~o Qu\^oc Bao is supported by the Vingroup Innovation Foundation under grant number VINIF.2023.TS.011. 
The work of Ph\`ung H\^o Hai is supported by the Vietnam Academy of Science and
Technology under grant number NVCC01.01/24-25. 
The work of D\`ao Van Thinh is supported by the Vingroup Innovation Foundation under grant number VINIF.2023.STS.49. }
\begin{document}
\maketitle
\begin{abstract}
Let $X$ be a smooth projective curve over a field $k$ of characteristic zero. 
The differential fundamental group of $X$ is defined as the Tannakian dual
to the category of vector bundles with (integrable) connections on $X$. This work investigates the relationship between the de Rham  cohomology of a vector bundle with connection and the group cohomology of the corresponding representation of the differential fundamental group of $X$. Consequently, we obtain some vanishing and non-vanishing results for the group cohomology. 
\end{abstract}


\section{Introduction}
For a smooth, geometrically connected algebraic variety $X$ over a field $k$ of characteristic $0$ and a rational point
$x \in X(k)$,  the differential fundamental group $\pi(X)$ of $X$ with base point in $x$ is an affine group scheme over $k,$ defined through the
Tannakian duality as follows. 
Consider the category   $\mathrm{MIC}^\mathrm{c}(X)$ of 
$\mathcal{O}_X$-coherent sheaves with $k$-linear integrable connections, i.e., an action
of the sheaf of differential operators $\mathcal D\textit{iff}(X)$. This is a $k$-linear abelian, 
rigid tensor category \cite{DM82}. The fiber functor at $x$:  
\begin{align*}
 x^*:   \mathrm{MIC}^\mathrm{c}(X) \longrightarrow \mathrm{Vec}_k,\quad 
    (\mathcal{V}, \nabla) \mapsto \mathcal{V}_{x}, 
\end{align*}
is faithfully exact. Hence,
Tannakian duality (cf. \cite[Theorem 2.11]{DM82}) applies and yields an affine group
scheme, denoted by $\pi(X)$, with the property that $x^*$ gives an equivalence:
\begin{align}\label{2}
 x^*:\mathrm{MIC}^\mathrm{c}(X) \overset{\cong}{\longrightarrow} 
 \mathrm{Rep}^\mathrm{f}(\pi(X)),
\end{align} 
where $\mathrm{Rep}^\textrm{f}(\pi(X))$ is the  category of finite dimensional
$k$-linear representations of $\pi(X)$.

If $k=\mathbb C$ the field of complex numbers, the Riemann-Hilbert correspondence
for projective curves tells us that $\pi(X)$ is the pro-algebraic hull of the topological
fundamental group of Riemann surfaces associated to $X$, which is well-known to be 
generated by $2g(X)$ elements subject to one relation ($g(X)$ is the genus of $X$). 

Recall that
the pro-algebraic completion $\widehat \Gamma_k$ of an abstract group $\Gamma$
is defined as the Tannakian dual of the category of finite dimensional 
complex linear representations of $\Gamma$. 
Unlike the \'etale fundamental group $\pi^\mathrm{et}(X)$ of $X$, 
$\pi(X)$ does not respect base
change and neither does the pro-algebraic group $\widehat\Gamma_k$. Therefore, it is
not known if $\pi(X)$ is isomorphic to the pro-algebraic completion of an abstract group in $2g(X)$
generators subject to one relation. Such an isomorphism, if it exists, is not compatible
with base change. A simple reason is that the exponential map which is needed to construct
the monodromy representations 
in the Riemann-Hilbert correspondence is defined only over $\mathbb C$ but not an arbitrary 
subfield of it. 

We aim to study $\pi(X)$ from a different perspective. We want to see
if $\pi(X)$ behave like an abstract group $\Gamma$ cohomologically. 
Our main result shows
that the group cohomology of $\pi(X)$ is much like that of the topological
fundamental group of the Riemann surface of the same genus. If there were a theory for the pro-algebraic closure of an abstract group parallel to Lubotsky-Magid \cite{LM82} we would 
have some interesting information on our $\pi(X)$. Unfortunately, such a theory
is not yet available. On the other hand, our result can be applied to compute
cohomology of the pro-unipotent quotient $\pi^\mathrm{uni}(X)$ of $\pi(X)$.
For this group, we can utilize Lubotsky-Magid's theory to conclude that   
the  $\pi^\mathrm{uni}(X)$ is a one-relator pro-unipotent group of rank $2g(X)$.

Let us briefly explain our main ideas. 
The equivalence in \eqref{2} induces a map from the
group cohomology $\mathrm{H}^{\bullet}(\pi(X),V)$ to 
the de Rham cohomology
$\mathrm{H}_\mathrm{dR}^{\bullet}(X, (\mathcal{V},\nabla))$ as follows.
First, the equivalence in \eqref{2}
extends to the categories of ind-objects (cf. \cite[\S~4]{Del89}). 
Thus, we have an equivalence
$$ x^*:\mathrm{MIC}^\mathrm{ind}(X) \overset{\cong}{\longrightarrow} 
 \mathrm{Rep}(\pi(X)),$$
where $\mathrm{MIC}^\mathrm{ind}(X)$ is the category of sheaves with integrable connections
which can be presented as inductive limits of objects from $\mathrm{MIC}^\mathrm{c}(X)$
(see also Remark \ref{rmk_01}).
 
An element $e$ of $\mathrm{H}^{n}(\pi(X),V)$ can be given in terms of an $n$-extension
of $k$ (as a trivial representation) 
by $V$ in $\mathrm{Rep}(\pi(X))$. Under the equivalence \eqref{2}, it yields
an $n$-extension of $(\mathcal O,\nabla)$ by $(\mathcal V,\nabla)$ in 
$\mathrm{MIC}^\mathrm{ind}(X)$. Then the connecting map 
$$\delta_e:\mathrm{H}_\mathrm{dR}^0(X, (\mathcal{O}, d))
\longrightarrow \mathrm{H}_\mathrm{dR}^n(X, (\mathcal{V}, \nabla))$$
yields an element $\delta_e(1)\in \mathrm{H}_\mathrm{dR}^n(X, (\mathcal{V}, \nabla))$. In this way,
we have defined a map
\begin{align}\label{22}
 \delta^{\bullet}(X): \mathrm{H}^\bullet(\pi(X),V) \longrightarrow \mathrm{H}_\mathrm{dR}^\bullet(X, (\mathcal{V}, \nabla)), \quad e\mapsto\delta_e(1).
\end{align}

These maps were studied in \cite{EH06}, where the authors showed
that for smooth, projective curves of genus at least 1, the maps $\delta^i, i\leq 2,$ are isomorphisms.  
Their proof of the surjectivity of $\delta^2$ (cf. \cite[Proposition 2.2]{EH06})
is incomplete, and the initial aim of this work is to provide a complete
proof of this claim.  
Our main finding is the following. 
\begin{thm}[Theorem \ref{theorem_delta_i}]
Let $X$ be a smooth projective geometrically connected curve with genus $g \geq 1.$ Then the maps $\delta^i$ are bijective for all $i \geq 0$ and all connections in 
$\mathrm{MIC}^\mathrm{ind}(X).$
\end{thm}
Since the de Rham cohomology on curves vanishes at degree $i\geq 3$, we obtain the same 
property for the cohomology of the differential fundamental groups of curves (of genus at 
least 1). This property is shared by the topological fundamental groups of Riemann
surfaces of genus at least 1. Our proof also mimics the method of proving that 
the fundamental group of a Riemann surface is ``algebraically good'' in \cite[4.3]{KPT09}. 
On the other hand, we notice that the claim of the theorem above fails trivially for 
zero genus curves because the
fundamental group is trivial in this case, while the second de Rham cohomology is non-vanishing.

Looking from the point of view of the algebraic $k(\pi,1)$ property (cf. \cite[6.2]{Sti13})
our finding says that smooth projective curves in characteristic $0$ are ``de Rham 
$k(\pi,1)$-spaces''. We notice that Di Proietto and Shiho
obtained similar results for ``log point'' \cite[section~5]{DS18}\footnote{We thank an anonymous referee for pointing these works out to us.}. 

Our method and result for $\pi(X)$ also apply 
to compute the cohomology of the pro-unipotent quotient $\pi^\mathrm{uni}(X)$ of $\pi(X)$. 
By results of Lubotsky and Magid \cite{LM82} we conclude that
$\pi^\mathrm{uni}(X)$ is a one-relator pro-unipotent group of rank $2g(X)$. 
For elliptic curves, we provide an explicit description of $\pi(X)$.

\noindent{\bf Notations}
\begin{itemize}
\item 
For an affine group scheme $G$ over $k$, let $\mathrm{Rep}(G)$
denote the category of $k$-linear representations.  
\item
The full subcategory of $\mathrm{Rep}(G)$ consisting of finite dimensional subrepresentations is denoted
by $\mathrm{Rep}^\mathrm{f}(G)$.

\item $\mathrm{MIC}(X)$ denotes the category of quasi-coherent $\mathcal O_X$-modules
with an integrable connection. The full subcategory of $\mathcal O_X$-coherent objects
is denoted by $\mathrm{MIC}^\mathrm{c}(X)$.
\item $\mathrm{MIC}^\mathrm{ind}(X)$ denotes the ind-category of $\mathrm{MIC}^\mathrm{c}(X)$
(i.e., the category of inductive systems of objects from $\mathrm{MIC}^\mathrm{c}(X)$, or equivalently, the category of left exact functors on the latter category to $\mathrm{Vec}_k$). 
\end{itemize}

\section{Integrable connections and de Rham cohomology}
Let $f: X \rightarrow \mathrm{Spec}(k)$ be a smooth geometrically connected 
variety where $k$ is a field of characteristic zero. 
 
\subsection{}
A connection $\nabla$ on a quasi-coherent $\mathcal O_X$-module 
$\mathcal V$ is a $k$-linear map
$$\nabla:\mathcal{D}\textit{er}(X)\longrightarrow 
\mathcal E\textit{nd}_k(\mathcal{V}).$$
We say that $\nabla$ is integrable  if
is a map of Lie algebras, or equivalently, if it extends to a map of sheaves of algebras
$$\nabla:\mathcal{D}\textit{iff}(X) \longrightarrow 
\mathcal E\textit{nd}_k(\mathcal{V}),$$
where $\mathcal{D}\textit{iff}(X)$ denotes  the sheaf of differential 
operators on $X$, cf. \cite[IV.16.8]{EGA4}.

We denote by $\mathrm{MIC}(X)$  the category of quasi-coherent $\mathcal{O}_X$-modules equipped
with $k$-linear integrable connections, and  
by $\mathrm{MIC}^\mathrm{c}(X)$ its full subcategory of $\mathrm{MIC}(X)$ 
consisting of coherent $\mathcal{O}_X$-modules with integrable connections.
Note that a coherent sheaf  equipped with an integrable connection is automatically
locally free (see, e.g., \cite{Kat70}).  

The structure sheaf $\mathcal O_X$ is equipped with a natural connection, denoted by $d$: 
the natural 
inclusion of $\mathcal{D}\textit{er}(X)$ into $\mathcal E\textit{nd}_k(\mathcal{O}_X)$.
It is a trivial connection. More generally, a module with 
connection is said to be trivial if it is the
direct sum of copies of $(\mathcal O_X,d)$. 

\subsection{}
The connection $\nabla$ yields a map, also called the connection and 
denoted by the same symbol,
$$\nabla:\mathcal V\longrightarrow \Omega_X^1\otimes_{\mathcal{O}_X}
\mathcal{V}. $$ 
For example, the map for the trivial connection on $\mathcal O_X$ is the differential
map
$$d:\mathcal O_X\longrightarrow \Omega_X^1.$$

By its integrability,  $\nabla$ extends to a differential on the complex 
$$\Omega^{\bullet}_{X}\otimes_{\mathcal{O}_X}\mathcal{V}.$$
We define the de Rham cohomology  $\mathrm{H}^{i}_\mathrm{dR}(X, (\mathcal{V},\nabla))$ of the module with connection $(\mathcal{V},\nabla)$ to be the hyper-cohomology:
$$ \mathrm{H}^i_\mathrm{dR}(X, (\mathcal{V},\nabla)):= 
\mathbf{H}^i(X, \Omega^{\bullet}_{X}\otimes_{\mathcal{O}_X}\mathcal{V}).$$
In particular, 
$$\mathrm{H}^0_\mathrm{dR}(X, (\mathcal{V}, \nabla)) = \mathrm H^0(X,\mathcal V^\nabla),$$ 
where  $\mathcal{V}^{\nabla}:=\mathrm{ker}\nabla$,   the kernel sheaf of 
$\nabla$, also called the
sheaf of horizontal sections of 
$(\mathcal{V},\nabla).$ 
According to \cite[(2.0)]{Kat70},  the de Rham
cohomology $\mathrm{H}^\bullet_\mathrm{dR}(X,(\mathcal{V},\nabla))$ can be 
identified  with the right derived functor of the  left exact functor:
\begin{align*}
 \mathrm{H}^0_\mathrm{dR}(X,-):   \mathrm{MIC}(X) &\longrightarrow \mathrm{Vec}_k.
\end{align*}
For each module with connection $(\mathcal{V},\nabla)$ we have a natural isomorphism 
$$\mathrm{Hom}_{\mathrm{MIC}(X)}((\mathcal{O}_X,d), (\mathcal{V},\nabla))
 \cong \mathrm H^0(X,\mathcal V^\nabla).$$
That is, we can identify the set of global horizontal sections of $(\mathcal V,\nabla)$ with the hom-set between the trivial module with connection $(\mathcal O_X,d)$ and  
 $(\mathcal V,\nabla)$. 
 
Since the category $\mathrm{MIC}(X)$ has enough injectives (cf. \cite[4.3.8]{ABC20}), 
the above isomorphisms induce isomorphisms of right derived functors 
\begin{align}\label{44}
 \mathrm{H}^{i}_\mathrm{dR}(X,-)\cong \mathrm{Ext}^{i}_{\mathrm{MIC}(X)}((\mathcal{O}_X,d),-).  
\end{align}
That is, an element in the $i$-th de Rham group of a module with 
connection $(\mathcal V,\nabla)$
can be given in terms of an $i$-extension of $(\mathcal O,d)$ by $(\mathcal V,\nabla)$.

\noindent{\bf Convention.}
When it is clear which connection on a vector bundle $\mathcal V$ we have in mind,
we shall omit it in notations, for example, we shall write $\mathrm{H}^i_\mathrm{dR}(X, \mathcal{V})$ instead of $\mathrm{H}^i_\mathrm{dR}(X, (\mathcal{V},\nabla))$.

\section{Cohomology of the differential fundamental group scheme
and the comparison map}

The category $\mathrm{MIC}^\mathrm{c}(X)$ is a rigid tensor $k$-linear abelian category, and when  a $k$-rational point $x\in X(k)$ is chosen,  we can equip 
$\mathrm{MIC}^\mathrm{c}(X)$ the fiber (i.e. faithfully exact) functor:
\begin{align*}
  x^*:  \mathrm{MIC}^\mathrm{c}(X) \longrightarrow \mathrm{Vec}_k,\quad
    (\mathcal{V}, \nabla) \mapsto \mathcal{V}|_{x},
\end{align*}
with values in the category of finite dimensional $k$-vector spaces. We notice that 
the faithfulness of $x^*$ follows from the assumption that $X$ is geometrically connected.

Tannakian duality \cite[Theorem 2.11]{DM82}  applies and yields an
affine group scheme $\pi(X)$ (we shall not mention
the base point $x$ in what follows) with the property that the fiber functor $x^*$ defines an equivalence of tensor categories,
\begin{equation}
\label{eq_tan1}
x^*:   \mathrm{MIC}^\mathrm{c}(X) \overset{ \cong}{\longrightarrow} \mathrm{Rep}^\mathrm{f}(\pi(X)).
\end{equation}
This equivalence extends to an equivalence between ind-categories (i.e., between the category of inductive direct limits of vector bundles with connection and the category of all
representations of $\pi(X)$):
\begin{equation}
\label{eq_tan2}
x^*:\mathrm{MIC}^\mathrm{ind}(X) \overset{ \cong}{\longrightarrow} \mathrm{Rep}(\pi(X)).
\end{equation}

\begin{rem}\label{rmk_01}
The ind-category $\mathrm{MIC}^\mathrm{ind}(X)$ of 
$\mathrm{MIC}^\mathrm{c}(X)$ can be identified with the category
of connections which are union (i.e. inductive limits)
of its coherent subconnections. Notice that
$\mathrm{MIC}^\mathrm{ind}(X)$ is strictly smaller than $\mathrm{MIC}(X)$.
For instance the sheaf of differential operators on $X$, 
$\mathcal{D}\textit{iff}(X)$ is a quasi-coherent module with integrable
connection (by the left action of itself) which is not the union of its 
coherent submodules with connections. In fact, a coherent submodule with 
connection cannot contain
the unit section of $\mathcal{D}\textit{iff}(X)$ as this section generates the whole 
$\mathcal{D}\textit{iff}(X)$. 
\end{rem}

\begin{rem}\label{rem-abelian}
If $g=1$, by choosing a $k$-point of $X$ to be the neutral element, $X$ is an elliptic curve and hence has a group structure. 
By K\"unneth theorem for $\pi(X)$ \cite[Corollary 10.47]{Del89}, it possesses another multiplication induced from the group structure on $X$. 
Using the Eckmann-Hilton argument, cf. \cite{EH62}, we conclude that $\pi(X)$ is commutative. 
\end{rem}

Let $V$ be an object of $ \mathrm{Rep}(\pi(X)).$ 
The invariant subspace of $V$ is defined to be ($k$ being equipped with
the trivial action):
\begin{align*}V^{\pi(X)}:=\mathrm{Hom}_{\pi(X)}(k,V).
\end{align*}
The fixed point functor $V\longmapsto (V)^{\pi(X)}$ is left exact functor (see \cite[2.10 -- (4)]{Jan87}). 
One defines the group cohomology $\mathrm{H}^{i}(\pi(X), V)$ as the $i$-th right derived functor of this functor.  Since $\mathrm{Rep}(\pi(X))$ has enough 
injectives, we get isomorphisms of right derived functors,  as stated in \cite[4.2]{Jan87}:
\begin{align}\label{33}
  \mathrm{H}^i(\pi(X),-) \cong \mathrm{Ext}^i_{\mathrm{Rep}(\pi(X))}(k, -).  
\end{align}


Let $(\mathcal{V},\nabla)$ be an object in  $\mathrm{MIC}^\mathrm{ind}(X),$ and let 
$V:=\mathcal V|_x$ be its fiber at $x$. The map $\delta^i$ mentioned in
\eqref{22} can be given as follows. An element  $e$ of $\mathrm{H}^i(\pi(X),V)$
is given in terms of an $i$-extension of $k$ by $V$, by means of the isomorphism 
in \eqref{33}. Under the equivalence \eqref{eq_tan2} it corresponds to an $i$-extension
of $(\mathcal O,d)$ by $(\mathcal V,\nabla)$ in $\mathrm{MIC}^\mathrm{ind}(X)$ and
thus can be considered as an extension in $\mathrm{MIC}(X)$, hence it
determines an element $\delta^i(e)$ of $\mathrm{H}^i_\mathrm{dR}(X,\mathcal{V})$.  

As $\mathrm{MIC}(X)$ is strictly larger than $\mathrm{MIC}^\mathrm{ind}(X)$, the
maps $\delta^i$ are in principle neither injective nor surjective. On the other hand,
since   $\mathrm{MIC}^\mathrm{c}(X)$ is a full subcategory of $\mathrm{MIC}(X)$
and is closed under extensions, $\delta^0$ and $\delta^1$ are
isomorphisms. The purpose of our work is to show that the maps $\delta^i$ for
$i\geq 2$ are also isomorphisms (in particular, $\delta^i=0$ for $i\geq 3$). 




\section{The main result}

We aim to show, in the case where $X$ is a projective smooth curve of genus $g\geq 1$ over $k$, that the maps $\delta^i$ are bijective for all $i\geq 2$. In particular, for
$i\geq 3$, $\delta^i$ is the zero map, that is
$$ \mathrm{H}^i(\pi(X),V)=0$$
for all  representations $V$. 

\begin{thm}\label{theorem_delta_i}
Let $X$ be a smooth projective geometrically connected curve with genus $g \geq 1$, equipped with a point $x\in X(k)$. Let $\pi(X):=\pi(X,x)$ be the
differential fundamental group at $x$. 
Then for any object  $(\mathcal V,\nabla)$ in $\mathrm{MIC}^\mathrm{ind}(X)$
and $V:=\mathcal V|_x$,  
the map
$$\delta^i:\mathrm{H}^i(\pi(X), V) \rightarrow 
 \mathrm{H}^i_\mathrm{dR}(X, (\mathcal V, \nabla)),$$ 
is bijective for all $i \geq 0.$
\end{thm}

For the proof, we have to distinguish between two cases, $g=1$ and $g\geq 2$.
We shall first need some lemmas.

\begin{lem}\label{lem-H2K}
    Let $X$ be a smooth projective geometrically connected curve with genus $g=1,$ equipped with a point $x \in X(k).$ Then 
    $$\mathrm{H}^2(\pi(X), k) \cong k.$$
Moreover, let $V$ be any representation of $\pi(X)$. Then 
$$ \mathrm{H}^3(\pi(X), V)=0.$$
\end{lem}
\begin{proof}
Since $X$ is an elliptic curve, the differential fundamental group $\pi(X)$ is commutative (see Remark \ref{rem-abelian}). Then we have
$$  \pi(X) = \pi^\mathrm{uni}(X) \times \pi^\mathrm{diag}(X),   $$
where $ \pi^\mathrm{uni}(X)$ is the pro-unipotent quotient of $\pi(X)$ and
$ \pi^\mathrm{diag}(X)$ is the pro-diagonal quotient of $\pi(X)$ (see, e.g., 
\cite[9.3]{Wat79}). The affine group scheme $\pi^\mathrm{uni}(X)$ corresponds
to the full subcategory of unipotent objects (i.e. iterated extensions of the
trivial module with connection $(\mathcal O_X,d)$), and the affine group scheme $\pi^\mathrm{diag}(X)$ corresponds to the semi-simple (diagonalizable)
connections (those which decompose into a direct sum of rank one connections after a base
change of the base field $k$).

Since the higher cohomology groups of $\pi^\mathrm{diag}(X)$ vanish \cite[Lemma 4.3-(b)]{Jan87},
using the Lyndon-Hochschild-Serre spectral sequence \cite[Proposition 6.6-(c)]{Jan87} associated to the exact sequence
$$1\longrightarrow \pi^\mathrm{uni}(X)\longrightarrow \pi(X)\longrightarrow
\pi^\mathrm{diag}(X)\longrightarrow 1,$$
we deduce that
\begin{align}\label{125}
\mathrm{H}^i(\pi^\mathrm{uni}(X),k) \cong \mathrm{H}^i(\pi(X), k), 
\end{align}
for all $i\geq 0$. 

By Tannakian duality,  $\delta^1$ is an isomorphism, hence
$$\mathrm{H}^1(\pi(X), k) \cong \mathrm{H}^1_\mathrm{dR}(X, \mathcal{O}_X) = k^2.$$
Since $\mathbb G_a$ has no non-trivial forms, cf.  \cite[14.58]{Mil17}, it follows from \cite[(1.2), (1.16)]{LM82} that $\pi^\mathrm{uni}(X) = \mathbb{G}_a \times \mathbb{G}_a.$ Therefore,
    $$\mathrm{H}^2(\pi(X), k) \cong k.$$
and 
$$ \mathrm{H}^3(\pi(X), V)=0,$$
for any $V$ in $\mathrm{Rep}(\pi(X)),$ cf. \cite[Remarks-(2), p.71]{Jan87}.
\end{proof}

\begin{lem}\label{lemh1non}
Let $X$ be a smooth projective curve of genus $g\geq 2.$ Let $(\mathcal V,\nabla)$ be a connection in $ \mathrm{MIC}^\mathrm{c}(X)$, $\mathcal V\neq 0$. 
Then the first de Rham cohomology $\mathrm{H}^1_\mathrm{dR}(X,(\mathcal{V},\nabla))$ is non zero. 
\end{lem}
\begin{proof} 
Let $r$ be the rank of $\mathcal V$.
We consider the first page of Hodge to de Rham spectral sequence:
 \begin{center}
 \begin{tikzcd}
\mathrm{H}^1(X,\mathcal{V}) \arrow{r}{\mathrm{H}^1(\nabla)} & \mathrm{H}^1(X,\Omega^1_{X/k}\otimes \mathcal{V}) \arrow{r} &0\\
\mathrm{H}^0(X,\mathcal{V}) \arrow{r}{\mathrm{H}^0(\nabla)} & \mathrm{H}^0(X,\Omega^1_{X/k} \otimes \mathcal{V}) \arrow{r} &0.
 \end{tikzcd}
\end{center}
According to \cite{Bro94}[(2.7)-p.164], the Euler characteristic of de Rham cohomology is equal to the Euler characteristic of the first page. That is, we have
$$ \chi_{dR}(\mathcal V) = \mathrm h^0(X,\mathcal{V}) -  \mathrm h^1(X,\mathcal{V}) +
\mathrm h^1(X,\Omega_{X/k}^1\otimes \mathcal{V}) -
\mathrm h^0(X,\Omega_{X/k}^1\otimes \mathcal{V}).  $$
Combine it with Serre's duality and using the Riemann-Roch theorem, we obtain
\begin{align*}
\chi_{dR}(\mathcal V) &=\mathrm h^0(X,\mathcal{V}) -  \mathrm h^1(X,\mathcal{V}) +
\mathrm h^0(X,\mathcal{V}^\vee) -\mathrm h^1(X, \mathcal{V}^\vee)\\
&=\mathrm{deg}(\mathcal{V}) + \mathrm{deg}(\mathcal{V}^{\vee}) + 2(1-g)r \\
& = 2(1-g)r<0. 
\end{align*}
%
%
Hence 
$$\mathrm h_\mathrm{dR}^1(X, \mathcal{V})=
\mathrm h_\mathrm{dR}^0(X, \mathcal{V}) -\chi_\mathrm{dR}(\mathcal{V}) > 0.$$
That is, $\mathrm{H}^1_\mathrm{dR}(X,(\mathcal{V},\nabla))\neq 0$. \end{proof}

\begin{lem}\label{lemkey}
Let $X$ be a smooth projective curve with genus $g\geq2,$ and $(\mathcal{V},\nabla)$ be an object in $\mathrm{MIC}^\mathrm{c}(X).$ Then there exists a connection $(\mathcal{V}',\nabla') \in \mathrm{Obj}(\mathrm{MIC}^\mathrm{c}(X))$ and an injective map $j: (\mathcal{V},\nabla) \hookrightarrow (\mathcal{V}',\nabla')$ such that 
$$\mathrm{H}^2_\mathrm{dR}(X,\mathcal{V}')=0.$$  
\end{lem}
\begin{proof} 
By Poincar\'e duality, it is equivalent to showing that there exists a surjective map
$\mathcal E\twoheadrightarrow\mathcal V^\vee$ of coherent connections such that  
$\mathrm{H}^{0}_\mathrm{dR}(X, \mathcal{E}) =0.$
By induction on the dimension of $\mathrm H^0_\mathrm{dR}(X,\mathcal E),$ it suffices to
show there exists  $\mathcal E$ surjecting on $\mathcal{V}^\vee$ and satisfying the strict inequality:
$$\mathrm{h}^{0}_\mathrm{dR}(X, \mathcal{E}) < 
\mathrm{h}^0_\mathrm{dR}(X,\mathcal{V}^\vee),$$
as long as $\mathrm{h}^0_\mathrm{dR}(X,\mathcal{V}^\vee)\neq 0$.

Thus, assume there exists an inclusion 
$(\mathcal O,d)\hookrightarrow (\mathcal V^\vee,\nabla)$ and let 
 $\mathcal{F}:= \mathcal V^{\vee}/\mathcal O $ be equipped with the quotient connection. 
 Consider an arbitrary, non-trivial,
simple connection 
($\mathcal{S}, \nabla_{\mathcal{S}}$) 
and apply the functor 
$\mathrm{Ext}^i_{\mathrm{MIC}(X)}(-,\mathcal{S})$
to the short exact sequence in $\mathrm{MIC}^\mathrm{c}(X):$
$$0 \longrightarrow \mathcal O \overset{e}{\longrightarrow} 
\mathcal V^{\vee} \longrightarrow \mathcal{F} \longrightarrow 0, $$
we obtain a long exact sequence:
$$  \mathrm{Ext}^1_{\mathrm{MIC}(X)}(\mathcal{F},\mathcal{S}) \longrightarrow 
\mathrm{Ext}^1_{\mathrm{MIC}(X)}(\mathcal{V}^{\vee},\mathcal{S})
\overset{\mathrm{ev}_{e}}{\longrightarrow} 
\mathrm{Ext}^1_{\mathrm{MIC}(X)}(\mathcal O,\mathcal{S})
\longrightarrow \mathrm{Ext}^2_{\mathrm{MIC}(X)}(\mathcal{F},\mathcal{S}).$$
By Poincar\'e duality we have 
\begin{align*}\mathrm{Ext}^2_{\mathrm{MIC}(X)}(\mathcal{F},\mathcal{S}) &\cong  \mathrm{H}^2_\mathrm{dR}(\mathcal O_X, \mathcal F^\vee\otimes
\mathcal S)\\
&\cong \mathrm{H}^0_\mathrm{dR}(\mathcal O_X,\mathcal F\otimes \mathcal S^\vee)^\vee\\
&\cong
\mathrm{Hom}_{\mathrm{MIC}^\mathrm{c}(X)}(\mathcal{S},\mathcal{F})^\vee.
\end{align*}

Since genus of $X$ is at least $1,$ there are infinitely many non isomorphic simple connections on 
$X$, for example, consider the trivial bundle equipped with a connection 
$\nabla:\mathcal O_X\longrightarrow \Omega_X,\quad 1\longmapsto \omega,$
for any non-zero form $\omega$. 
Since $\mathrm{MIC}^\mathrm{c}(X)$ is equivalent to 
$\mathrm{Rep}^\mathrm{f}(\pi(X))$, we can choose a non-trivial simple connection 
$(\mathcal{S}, \nabla_{\mathcal{S}})$  such that 
$$\mathrm{Hom}_{\mathrm{MIC}^\mathrm{c}(X)}(\mathcal{S},\mathcal{F})^\vee= 0,$$ 
consequently, the map 
$$\mathrm{ev}_e:\mathrm{Ext}^1_{\mathrm{MIC}(X)}
(\mathcal{V}^{\vee},\mathcal{S}) {\longrightarrow} 
\mathrm{Ext}^1_{\mathrm{MIC}(X)}(\mathcal O,\mathcal{S})$$
is surjective. 
 
According to Lemma \ref{lemh1non}, there exists  
a non-split extension of connections:
$$\varepsilon: \mathcal S\longrightarrow \mathcal G\longrightarrow \mathcal O.$$
With the assumption that $(\mathcal S,\nabla_{\mathcal S})$ 
is non-trivial and simple, we have  
that $\mathrm{H}^0_\mathrm{dR}(X, \mathcal G)=0$.   
Let $\epsilon$ be a preimage of $\varepsilon$ along the surjective map   
$\mathrm{ev}_{e}$.  
In terms of extensions, these are related by the following diagrams in $\mathrm{MIC}^\mathrm{c}(X)$
$$\xymatrix{ 
\qquad\qquad \epsilon: 0\ar[r]& \mathcal{S}\ar[r] \ar@{=}[d]& 
 \mathcal{E}\ar[r]& \mathcal{V}^{\vee} 
 \ar[r]& 0\\ 
 \varepsilon= \mathrm{ev}_e(\epsilon): 0\ar[r]&  \mathcal{S}  \ar[r]&
 \mathcal G\ar[r]\ar[u]& \mathcal O\ar[r]\ar[u]_e&0.}$$
Take the associated long exact sequence of de Rham cohomology we have
$$\xymatrix{
0=\mathrm{H}^{0}_\mathrm{dR}(X,\mathcal S) \ar[r]&
\mathrm{H}^{0}_\mathrm{dR}(X,\mathcal E)\ar[r] &
\mathrm{H}^{0}_\mathrm{dR}(X,\mathcal V^\vee)\ar[r]&
\mathrm{H}^{1}_\mathrm{dR}(X,\mathcal S)\\
0=\mathrm{H}^{0}_\mathrm{dR}(X,\mathcal S) \ar[r]\ar@{=}[u]&
\mathrm{H}^{0}_\mathrm{dR}(X,\mathcal G) =0\ar[r]\ar[u] &
\mathrm{H}^{0}_\mathrm{dR}(X,\mathcal O)\ar@{^(->}[r]\ar@{^(->}[u]&
\mathrm{H}^{1}_\mathrm{dR}(X,\mathcal S).\ar@{=}[u]
}$$
Since $\mathrm{H}^{0}_\mathrm{dR}(X,\mathcal O)=k\neq 0$, the rightmost
upper horizontal map
$$\mathrm{H}^{0}_\mathrm{dR}(X,\mathcal V^\vee)\longrightarrow
\mathrm{H}^{1}_\mathrm{dR}(X,\mathcal S)$$
should be non-zero, that is
$\mathrm{H}^{0}_\mathrm{dR}(X,\mathcal E)\neq 
\mathrm{H}^{0}_\mathrm{dR}(X,\mathcal V^\vee),$
whence the desired inequality. \end{proof}

\begin{proof}[Proof of Theorem \ref{theorem_delta_i}]
By Tannakian duality, $\delta^0$ and $\delta^1$ are isomorphisms for
connections $(\mathcal V,\nabla)$ in $\mathrm{MIC}^\mathrm{ind}(X)$.  
Therefore, by the ``dimension shifting'' trick, we conclude that $\delta^2$ is injective.
This has already been done in \cite{EH06}, we recall it here for  completeness.

\textit{Injectivity of $\delta^2$.}
Let $(\mathcal J,\nabla)$ be an injective envelope of $(\mathcal V,\nabla)$ in
$\mathrm{MIC}^\mathrm{ind}(X)$ and denote by $J$ its fiber at $x$.
Thus, $J$ is an injective envelope of $V$ in $\mathrm{Rep}(\pi(X)).$ 
Since both de Rham cohomology and group cohomology commute with direct limits, we have 
$$
0=\mathrm{H}^1(\pi(X),J)\cong \mathrm{H}^1_\mathrm{dR}(X,\mathcal{J}).$$ 
The long exact sequences associated to the exact sequences
$$0\longrightarrow \mathcal V\longrightarrow \mathcal J\longrightarrow
\mathcal{J/V}\longrightarrow 0$$ and 
$$0\longrightarrow V\longrightarrow J\longrightarrow J/V\longrightarrow 0$$
yield
$$\xymatrix{
0=\mathrm{H}^1(\pi(X),J)\ar[r]\ar[d]^\cong& \mathrm{H}^1(\pi(X),J/V)\ar[r]\ar[d]^\cong
& \mathrm{H}^2(\pi(X),V)\ar[r]\ar@{^(->}[d]^{\delta^2}& 0\\
\mathrm{H}^1_\mathrm{dR}(X,\mathcal{J})\ar[r]&
\mathrm{H}^1_\mathrm{dR}(X,\mathcal{J/V})\ar[r]& 
\mathrm{H}^2_\mathrm{dR}(X,\mathcal{V})\ar[r]&\ldots }$$

\textit{Surjectivity of $\delta^2$.} There are two cases: the genus is 1
or larger than 1. 

Case $g\geq 2$: 
Let $(\mathcal{V},\nabla)$ be a connection in $\mathrm{MIC}^\mathrm{c}(X)$ and let $(\mathcal{V},\nabla) \hookrightarrow (\mathcal{V}',\nabla')$ as in Lemma \ref{lemkey}, that is,  $   {\mathrm{H}^2_\mathrm{dR}(X,\mathcal{V}')}=0$.
The long exact sequences of cohomology associated to the exact sequences
$0 \rightarrow \mathcal{V}  \rightarrow 
\mathcal{V}' \rightarrow \mathcal{V}'/\mathcal{V} \rightarrow 0,$ in $\mathrm{MIC}(X)$
and of their fibers in $\mathrm{Rep}^\mathrm{f}(\pi(X)):$
 $0 \rightarrow V \rightarrow V' \rightarrow V'/V \rightarrow 0  $  fit in the 
following diagram:
$$\xymatrix{
	{\mathrm{H}^1(\pi(X),V'/V)}\ar[r]\ar[d]^{\delta^1}_\cong & {\mathrm{H}^2(\pi(X),V)} \ar[r]\ar@{^(->}[d]^{\delta^2}& {\mathrm{H}^2(\pi(X),V')} \ar@{^(->}[d]^{\delta^2}
\\
	{\mathrm{H}^1_\mathrm{dR}(X,\mathcal{V}'/\mathcal{V})} \ar@{->>}[r]& {\mathrm{H}^2_\mathrm{dR}(X,\mathcal{V})} \ar[r]& {\mathrm{H}^2_\mathrm{dR}(X,\mathcal{V}')}=0. }$$
We conclude that the middle arrow is surjective.  


Case  $g=1$. We prove by induction on the rank of the vector bundle. Notice that the claim holds for the trivial connection $(\mathcal O_X,d)$.
Indeed, according to Lemma \ref{lem-H2K} 
$\mathrm{H}^2(\pi(X), k)=k=\mathrm{H}^2_\mathrm{dR}(X,\mathcal O_X)$.
Hence $\delta^2$, being injective, is bijective at $(\mathcal O_X,d)$. 

Consider a vector bundle with connection $\mathcal V$.
If $\mathrm{H}^2_\mathrm{dR}(X,\mathcal{V})= 0$, there is nothing to prove. 
Thus assume that $\mathrm{H}^2_\mathrm{dR}(X,\mathcal{V})\neq 0$.
By Poincar\'e duality,
$\mathrm{H}^0_\mathrm{dR}(X,\mathcal{V}^\vee)\neq 0$. By duality, we have
an exact sequence in $\mathrm{MIC}(X):$
$$0 \rightarrow \mathcal{G} \rightarrow \mathcal{V} \rightarrow \mathcal{O}_X \rightarrow 0.$$
This yields the corresponding exact sequence on their fibers as representations of 
$\pi(X)$ and hence the long exact sequences, fitting in the following diagram:
$$\xymatrix{
\mathrm{H}^1(\pi(X), k)\ar[r]\ar[d]^\cong&
\mathrm{H}^2(\pi(X), G)\ar[r]\ar[d]^\cong&
\mathrm{H}^2(\pi(X), V)\ar[r]\ar[d]&
\mathrm{H}^2(\pi(X), k)\ar[r]\ar[d]^\cong& \mathrm{H}^3(\pi(X), G)\ar@{=}[d]\\
\mathrm{H}^1_\mathrm{dR}(X,\mathcal O_X)\ar[r]&
\mathrm{H}^2_\mathrm{dR}(X,\mathcal G)\ar[r] &
\mathrm{H}^2_\mathrm{dR}(X,\mathcal V)\ar[r]&
\mathrm{H}^2_\mathrm{dR}(X,\mathcal O_X)\ar[r]&
0.}$$
According to Lemma \ref{lem-H2K}, $\mathrm{H}^3(\pi(X), G)=0$. 
Thus, if the second vertical map is bijective, so is the third vertical map. This finishes the
induction step, hence finishes the proof that $\delta^2$ is an isomorphism at objects of $\mathrm{MIC}^\mathrm{c}(X)$.

As de Rham cohomology and group cohomology commute with direct limits, we see that $\delta^2$
is an isomorphism at objects of $\mathrm{MIC}^\mathrm{ind}(X).$

The last step is to show that $\delta^i$ are isomorphisms for all $i\geq 0$.
Since the de Rham cohomology vanishes in 
degrees $i\geq 3$, it suffices to show that $\delta^i$ is injective. 
We employ the ``dimension shiting'' trick
 -- consider an injective envelope $(J,\nabla)$ of $(\mathcal V,\nabla)$
in $\mathrm{MIC}^\mathrm{ind}(X)\cong\mathrm{Rep}(\pi(X))$
and the associated long exact sequences of de Rham cohomology and group 
cohomology. Then, we conclude by
induction that $\delta^i$ ($i\geq 3$) is injective:
$\mathrm{H}^i(\pi(X),V)\cong \mathrm{H}^{i-1}(\pi(X),J/V)
\cong \mathrm{H}^{i-1}_\mathrm{dR}(X,\mathcal{J/V})\hookrightarrow 
\mathrm{H}^{i}_\mathrm{dR}(X,\mathcal{V}).$
\end{proof} 

Since the de Rham cohomology of curves vanishes in any degree larger than $2$ we
obtain immediately: 
\begin{cor}
For a smooth curve $X$ of genus at least 1, the group 
cohomology of $\pi(X)$ with coefficients in any representation vanishes in degree
$n>2$.
\end{cor}

\section{Cohomology of the unipotent group}
The pro-unipotent quotient $\pi^\mathrm{uni}(X)$
of $\pi(X)$ is defined to be the Tannakian dual to
the category of nilpotent objects, that is, vector bundles with connections that are iterated 
extension of the trivial connection $(\mathcal O_X,d)$. 
According to \cite[Corollary~10.43]{Del89} $\pi^\mathrm{uni}(X)$ is compatible with
base change. 
As an application of our main result, we give a purely algebraic computation of the
cohomology of $\pi^\mathrm{uni}(X)$.  
Then a beautiful result of Lubotsky-Magid \cite{LM82} yields the structure
of $\pi^\mathrm{uni}(X)$.

Recall that a pro-unipotent group $G$ is said to be finitely generated 
if there is a finite set of elements of $G$, the abstract subgroup generated by which is
Zariski dense in $G$. The minimal cardinal of such a set is called the rank of $G$.
This number is equal to $r_1:=\dim_k\mathrm{H}^1(G,k)$.
In this case, there is an exact sequence $1\longrightarrow R\longrightarrow F\longrightarrow G\longrightarrow1$ of pro-unipotent groups in which $F$ is free of the same rank as $G$.
$R$ is also a free pro-unipotent group. $R$ is generally infinitely generated.
If $\mathrm{H}^2(G,k)$ is finite dimensional
then $R$ is finitely generated \textit{as a normal subgroup of rank} $r_2:=\dim_k\mathrm{H}^2(G,k)$. The latter amounts to saying that there are $r_2$ elements in $R$ whose conjugates
by $F$ generate $R$, further $r_2$ is the minimum number of such elements. 
For more detail the reader is referred to \cite{LM82}. 
\begin{prop}
Let $X$ be a smooth projective curve of genus at least $1.$ Then for any
representation $V$ of $\pi^\mathrm{uni}(X),$ the canonical maps
$$\gamma^i_V: \mathrm{H}^i(\pi^\mathrm{uni}(X),V)\longrightarrow
\mathrm{H}^i(\pi(X),V), i=0,1,\ldots$$
are bijective. In particular $\mathrm{H}^i(\pi^\mathrm{uni}(X),V)$ 
vanishes for all $V$ and all $i\geq 3$. Consequently, when $k$ is algebraically
closed, $\pi^\mathrm{uni}(X)$ is
a one-relator pro-unipotent group of rank  $2g(X)$. 
\end{prop}
\begin{proof}

It is obvious that $\gamma^0_V$ and $\gamma^1_V$ are bijective.  Then, using
the ``dimension shifting'' trick as at the end of the proof of Theorem
\ref{theorem_delta_i}, we see that $\gamma^2$ is injective at all nilpotent representations.

Next, we claim that $\gamma^2$ is surjective at the trivial representation. Indeed,
employing the Poincar\'e duality for de Rham cohomology and Theorem
\ref{theorem_delta_i}, we can find two extensions of the trivial connection 
$(\mathcal O,d)$ such that their cup product is a non-trivial two-extension.
That means, in the following commutative diagram:
$$\xymatrix{
\mathrm{H}^1(\pi^\mathrm{uni}(X),k)\times \mathrm{H}^1(\pi^\mathrm{uni}(X),k)
\ar[rr]^{\qquad\cup}
\ar[d]^\cong &&
\mathrm{H}^2(\pi^\mathrm{uni}(X),k)\ar[d]\\
\mathrm{H}^1(\pi(X),k)\times \mathrm{H}^1(\pi(X),k)\ar@{->>}[rr]^{\qquad\cup} &&
\mathrm{H}^2(\pi(X),k)\cong k}$$
the lower horizontal map is surjective. As the left vertical map is bijective,
we conclude that the right vertical map is surjective. 
Thus we conclude that $\mathrm{H}^2(\pi^\mathrm{uni}(X),k)\cong k$.
According to \cite[Theorem~3.14]{LM82}, $\pi^\mathrm{uni}(X)$ has 
cohomological dimension 2, that is $\mathrm{H}^i(\pi^\mathrm{uni}(X),V)=0$
for any $V$ and any $i\geq 3$.
Recall that we also have $\mathrm{H}^i(\pi(X),V)=0$
for any $i\geq 3$.

Now, for any representation $V$ of $\pi^\mathrm{uni}(X)$, as $V$
is an iterated extension of the trivial representation, an induction argument together with
the long exact sequence of cohomology implies that
$$\mathrm{H}^2(\pi^\mathrm{uni}(X),V)\cong \mathrm{H}^2(\pi(X),V).$$
According to \cite[Theorem 3.2 and 3.11]{LM82}, when $k$ is algebraically closed,
$\pi^\mathrm{uni}(X)$ is  a  one-relator pro-unipotent group of rank $2g(X)$. 
\end{proof}

\begin{rem}We use the ``dimension shifting'' trick in the above proof 
to conclude the injectivity of the map $\delta^2$. But we cannot show that
$\delta^2$ is bijective as in the proof of Theorem \ref{theorem_delta_i}, 
hence we cannot go on to show that
$\mathrm{H}^i(\pi^\mathrm{uni}(X),V)=0$
for $i\geq 3$. Our proof is really based on the result of Lubotsky-Magid 
\cite[Theorem~3.14]{LM82},
which is a very specific property possessed by pro-unipotent group schemes.
\end{rem}

\begin{rem}
Assume that $g=1$ and 
the field $k$ is algebraically closed. Then $\pi(X)$ is commutative, hence decomposes into 
the product of the pro-diagonal
part $\pi^\mathrm{diag}(X)$ and the pro-unipotent part 
$\pi^\mathrm{uni}(X) = \mathbb{G}_a \times \mathbb{G}_a$ (cf. proof of Lemma
 \ref{lem-H2K}).
The pro-diagonal part $\pi^\mathrm{diag}(X)$ is 
Tannakian dual to the semi-simple category
generated by line bundles equipped with connections 
(automatically integrable as we are on a curve). Let $\mathbb X$ be the set of isomorphism class of
line bundles with connections, it is a commutative group with the multiplication induced from the
tensor product. Then $$\pi^\mathrm{diag}(X)\cong \mathrm{Diag}(\mathbb X):=\mathrm{Spec}(k[\mathbb X]),$$ 
where $k[\mathbb X]$ denotes the group ring of $\mathbb X$. It is well-known that
a line bundle on $X$ carries a connection if and only if it has degree 0 (cf. \cite[7.4]{Kat72}). More precisely, $\mathbb X$ is isomorphic to the
first hypercohomology group $\mathbb H^1(X,\Omega^*_X)$, where
$\Omega^*_X$ is the multiplicative de Rham complex, which for curves reads
$\mathcal O_X^*\xrightarrow{\mathrm{dlog}}\Omega^1_X.$
There is an exact sequence
$$0\longrightarrow \mathrm{H}^0(X,\Omega^1_X)\longrightarrow 
\mathbb H^1(X,\Omega^*_X)\longrightarrow \mathrm{Pic}^0(X)\longrightarrow 0,$$
where the last term represents the set of line bundles of degree 0 and the first term
represents the difference between choices of connections on a given line bundle. Notice that
$X\cong\mathrm{Pic}^0(X)$ as $X$ is an elliptic curve, the canonical map is given by
$P\longmapsto \mathcal O(P-O)$. The reader may compare with the 
differential fundamental group 
of the punctured affine line, where it is the product of the diagonal
group on the group $k/\mathbb Z$ and a copy of $\mathbb G_a$,
cf. \cite{HdST24}.  
\end{rem}

\begin{center}\bf 
Acknowledgments 
\end{center}

We thank H\'el\`ene Esnault for engaging in insightful discussions and Jo\~ao Pedro 
dos Santos for explaining to us about connections on 
line bundles on elliptic curves.
We thank the anonymous referee for his/her careful reading of the work and
many critics that greatly improved the presentation of our work.

\bibliographystyle{alpha}

\end{document}